\documentclass[letterpaper, 10 pt, conference]{ieeeconf}

\IEEEoverridecommandlockouts
\usepackage{amsmath, amssymb, amsthm}
\usepackage{graphicx}         
\usepackage{cite}
\usepackage{svg}
\usepackage{algorithm}
\usepackage{algpseudocode}
\algrenewcommand\algorithmicrequire{\textbf{Input:}}

\newtheorem{theorem}{Theorem}
\newtheorem{lemma}[theorem]{Lemma}

\theoremstyle{definition}

\newtheorem{assumption}{Assumption}
\newtheorem{problem}{Problem}
\theoremstyle{remark}
\newtheorem{remark}{Remark}
\theoremstyle{definition}

\newcommand{\R}{\mathbb{R}}
\newcommand{\N}{\mathbb{N}}

\usepackage{amsmath,amssymb,amsthm}
\newtheorem{prop}{Proposition}

\newtheorem{example}{Example}

\providecommand{\R}{\mathbb{R}}
\providecommand{\N}{\mathbb{N}}
\providecommand{\inner}[2]{\langle #1, #2 \rangle}
\providecommand{\SOL}{\operatorname{SOL}}
\providecommand{\dist}{\operatorname{dist}}

\providecommand{\Normal}{\mathrm{N}}

\providecommand{\E}{\mathbb{E}}

\title{\LARGE \bf
Tikhonov regularization for monotone v-GNE problems with non-linear constraints
}

\author{Tatiana Tatarenko$^{1}$ and Lucas Wey Hacker$^{2}$ %
\thanks{$^{1,2}$ T. Tatarenko and L.W. Hacker are with the Department of Control Methods and Intelligent Systems, Technical University of Darmstadt, Germany {\tt\small tatiana.tatarenko@tu-darmstadt.de, lucas.hacker@tu-darmstadt.de}. This work was funded by the Deutsche Forschungsgemeinschaft (DFG,
German Research Foundation), project number 528033031.}}%

\begin{document}

\maketitle
\thispagestyle{empty}
\pagestyle{empty}

\begin{abstract}
    Generalized Nash equilibrium problems model noncooperative agents whose decisions are coupled through shared constraints. This paper addresses the computation of variational equilibria of such games when the pseudo-gradient is merely monotone and the coupling constraints are convex but nonlinear. In this setting, the variational equilibrium does not need to be unique, and the projection onto the feasible set requires the solution of a constrained optimization problem at every iteration. We propose a randomized algorithm that combines a Tikhonov-regularized forward step with a projection onto the local action sets intersected with the halfspace obtained by linearizing a single, randomly selected coupling constraint. This projection is computed exactly by a one-dimensional search, and its cost does not depend on the number of coupling constraints; in particular, no projection onto the nonlinear feasible set is required. When the constraints satisfy an error bound (e.g., under Slater's condition or when they are affine), we show that the expected distance of the iterates from the feasible set is of the order of the step size. We further prove that the iterates converge in mean square to the minimum-norm variational equilibrium, without assuming strong monotonicity of the pseudo-gradient.
\end{abstract}
\section{Introduction}\label{sec:intro}

Generalized Nash equilibrium problems (GNEPs) model noncooperative agents whose feasible sets depend on the decisions of the others, typically through shared resources \cite{facchinei2007gnep,fischer2014gnep}, and are widely used in control applications such as resource allocation, congestion control and energy systems \cite{yi2019operator,belgioioso2021distributed}. Of particular interest are the variational equilibria (v-GNEs), i.e., the solutions of a single variational inequality (VI) over the shared feasible set \cite{facchinei2007vi}: at a v-GNE all players pay the same price for each shared constraint \cite{kulkarni2012variational}, as in Rosen's normalized equilibria \cite{rosen1965existence}. We address the computation of v-GNEs in a setting that combines two following difficulties.

The first one is the lack of uniqueness. Most v-GNE seeking algorithms require a strongly monotone pseudo-gradient \cite{yi2019operator, TatKamECC25}, which guarantees a unique v-GNE. Under mere monotonicity, which arises, e.g., when the costs are convex but not strictly convex in the players' own actions, the set of v-GNEs is convex but may be a continuum. Algorithms for this setting exist, e.g., the extragradient method \cite[Sec.~12.1.2]{FaccPang2}, but the equilibrium they reach depends on the initialization and on the algorithm parameters. Hence, the outcome is neither predictable nor reproducible, and different equilibria may induce very different costs. This calls for an equilibrium selection rule, i.e., for computing the v-GNE that minimizes a secondary criterion over the solution set of the VI. 
The simplest criterion is the norm: Tikhonov regularization, which perturbs the operator by a vanishing multiple of the identity, selects the minimum-norm solution \cite{BakPolyak74}, \cite[Sec.~12.2]{FaccPang2}, and perturbing it by a vanishing multiple of $x-x_{\mathrm{r}}$ instead selects the solution closest to a reference point $x_{\mathrm{r}}$, e.g., a nominal operating point. This selection is unique, independent of both the algorithm and its initialization, and requires no knowledge of the solution set.

The second difficulty concerns the coupling constraints. Projection methods, including the original Tikhonov regularization procedure \cite{BakPolyak74}, \cite[Sec.~12.2]{FaccPang2}, require the projection onto the feasible set at every iteration \cite[Ch.~12]{FaccPang2}, which, for nonlinear coupling constraints, is itself a constrained optimization problem. Primal-dual and operator-splitting methods avoid this projection by dualizing the coupling constraints, but they are typically developed for affine coupling constraints \cite{yi2019operator, TatKamECC25}. For nonlinear coupling constraints, the lifted operator is monotone but not Lipschitz continuous, so that line-search or adaptive step-size rules are required, e.g., the golden ratio method \cite{malitsky2020golden}. Since such step sizes are computed from global quantities of the iterates, they are difficult to distribute among the players. Moreover, dualization lifts the problem to the space of primal and dual variables, whose dimension grows with the number of coupling constraints, which may slow down convergence.

\emph{Contributions.} We address both difficulties simultaneously, staying in the primal space and without complex projections onto the feasible set.
\begin{enumerate}
  \item We propose a randomized algorithm for computing the minimum-norm v-GNE of merely monotone games with convex, possibly nonlinear and nonsmooth, coupling constraints (Algorithm~\ref{alg:main}). Each iteration takes a forward step on the Tikhonov-regularized pseudo-gradient and projects onto the local action sets intersected with the linearization of one randomly selected coupling constraint, in the spirit of random projection methods \cite{polyak2001random,nedic2011random}; its cost does not depend on the number of coupling constraints, and a fully explicit variant is given (Remark~\ref{rem:sequential}).
  \item Under an error bound on the constraint system, implied, e.g., by Slater's condition or by affine coupling constraints (Lemma~\ref{lem:errorbound_suff}), we show asympotic feasibility of the algorithm, i.e., that the expected distance of the iterates from the feasible set is of the order of the chosen diminishing step size (Lemma~\ref{lem:feasibility}).
  \item We prove mean-square convergence to the minimum-norm v-GNE (Theorem~\ref{thm:main}) for a pseudo-gradient that is merely continuous and monotone, i.e., without strong monotonicity.
\end{enumerate}
Since the game enters the analysis only through the pseudo-gradient, the results apply to monotone VIs over sets of the same form. However, the game structure is essential for implementation: the forward step can be performed by each player on its own action, and, in the explicit variant, the projections are componentwise up to a single scalar shared by all players. This makes the algorithm a natural candidate for distributed implementations, which we leave for future work.
\section{Preliminaries and problem formulation}\label{sec:prelim}

\emph{Notation.} We let $[s]_+ := \max\{s,0\}$ for $s\in\R$, and equip $\R^n$ with the inner product $\inner{x}{y} := x^\top y$ and the norm $\|x\| := \inner{x}{x}^{1/2}$; $I$ is the identity, $(x_1,\ldots,x_n)$ is a column vector, and vector inequalities are componentwise. For a nonempty, closed and convex set $C\subseteq\R^n$, $\Pi_C(x) := \operatorname{argmin}_{y\in C}\|y-x\|$ is the (single-valued) projection onto $C$ and $\dist(x,C) := \|x-\Pi_C(x)\|$. For all $y\in\R^n$ and $z\in C$ \cite[Ch.~3]{bauschke2017convex},
\begin{equation}\label{eq:fejer_basic}
  \|\Pi_C(y)-z\|^2\le\|y-z\|^2-\|y-\Pi_C(y)\|^2.
\end{equation}
A mapping $F:\R^n\to\R^n$ is monotone on $U\subseteq\R^n$ if $\inner{F(x)-F(y)}{x-y}\ge 0$ for all $x,y\in U$, and strongly monotone on $U$ with modulus $c>0$ if $\inner{F(x)-F(y)}{x-y}\ge c\|x-y\|^2$ for all $x,y\in U$. The variational inequality $\mathrm{VI}(C,F)$ consists in finding $x^\star\in C$ such that $\inner{F(x^\star)}{x-x^\star}\ge 0$ for all $x\in C$; its solution set is denoted by $\SOL(C,F)$.

\subsection{Generalized Nash equilibrium problems}\label{subsec:gnep}

We consider a game with players $\nu\in\mathcal{N} := \{1,\ldots,n\}$ \cite{facchinei2007gnep}. Player~$\nu$ chooses an action $x_\nu\in\R$ (scalar actions only simplify the notation; the analysis extends verbatim to $x_\nu\in\R^{d_\nu}$), we write $x := (x_1,\ldots,x_n)$ and $x=(x_\nu,x_{-\nu})$, and player~$\nu$ minimizes a cost $J_\nu(x_\nu,x_{-\nu})$, where $J_\nu(\cdot,x_{-\nu})$ is convex and differentiable for every $x_{-\nu}$. The \emph{local constraints} are $x_\nu\in A_\nu := [a_\nu,b_\nu]$, $a_\nu\le b_\nu$, so that $A := A_1\times\cdots\times A_n$ is a box and
\begin{equation}\label{eq:proj_A}
  \big(\Pi_A(x)\big)_\nu = \min\big\{b_\nu,\max\{a_\nu,x_\nu\}\big\},\quad \nu\in\mathcal{N}.
\end{equation}
The \emph{coupling shared constraints} are $g(x)\le 0$, where $g=(g_1,\ldots,g_m):\R^n\to\R^m$ and each $g_i$, $i\in\mathcal{M} := \{1,\ldots,m\}$, is convex,  hence continuous \cite[Cor.~10.1.1]{rockafellar1970convex}, and possibly \textit{nonlinear}. We let
\begin{equation}\label{eq:coupling_sets}
  X_i := \{x : g_i(x)\le 0\},\quad X := \bigcap_{i\in\mathcal{M}} X_i,\quad S := X\cap A.
\end{equation}
Given $x_{-\nu}$, player~$\nu$ solves $\min_{x_\nu} J_\nu(x_\nu,x_{-\nu})$ subject to $x_\nu\in K_\nu(x_{-\nu}) := \{x_\nu\in A_\nu : g(x_\nu,x_{-\nu})\le 0\}$. A point $x^\star$ is a \emph{generalized Nash equilibrium} (GNE) if, for every $\nu\in\mathcal{N}$, $x^\star_\nu$ solves this problem for $x_{-\nu}=x^\star_{-\nu}$. The \emph{pseudo-gradient} of the game is
\begin{equation}\label{eq:pseudo_gradient}
  F(x) := \big( \partial J_1(x)/\partial x_1,\ \ldots,\ \partial J_n(x)/\partial x_n \big).
\end{equation}

The following assumptions hold throughout the paper.

\begin{assumption}[Standing assumptions]\label{assum:standing}\leavevmode
\begin{enumerate}
  \item[(A1)] $F$ is monotone on $A$;
  \item[(A2)] $A$ is compact and $X$ is nonempty, so that $S$ is nonempty and compact;
  \item[(A3)] $F$ is continuous on $A$;
  \item[(A4)] each $g_i$ is differentiable on $\R^n$;
  \item[(A5)] (error bound) for every bounded $B\subset\R^n$, there is $\tau_B\ge 1$ such that, for all $x\in B$,
  \begin{equation}\label{eq:errorbound}
    \dist(x,S)\le\tau_B\max\Big\{\max_{i\in\mathcal{M}}[g_i(x)]_+,\ \dist(x,A)\Big\}.
  \end{equation}
\end{enumerate}
\end{assumption}

Monotonicity is imposed on $A$ rather than on $S$, because the iterates of our algorithm belong to $A$ but, in general, not to $S$. By (A4) and \cite[Cor.~25.5.1]{rockafellar1970convex}, each $\nabla g_i$ is continuous; (A4) can be dispensed with Remark~\ref{rem:nonsmooth}. (A5) is a regularity condition on the constraint system $\{g(x)\le 0,\ x\in A\}$: it quantifies how small violations of the constraints translate into proximity to $S$. It holds, e.g., under Slater's condition or for affine $g$, as we show next.

We say that \emph{Slater's condition} holds if
\begin{equation}\label{eq:slater}
  \text{there is } \bar x\in A \text{ with } g_i(\bar x)<0 \text{ for all } i\in\mathcal{M}.
\end{equation}
In this case, we let $\bar\delta := -\max_{i\in\mathcal{M}} g_i(\bar{x}) > 0$ be the Slater margin and, for a bounded set $B\subset\R^n$, $\rho_B := \sup_{x\in B}\|x-\bar{x}\|$.

\begin{lemma}\label{lem:errorbound_suff}
Let (A2) hold. Then (A5) holds in each of the following cases.
\begin{enumerate}
  \item\label{lem:eb_slater} Slater's condition~\eqref{eq:slater} holds. In this case, for every nonempty bounded $B\subset\R^n$, there is $\kappa_{B,A}\ge 1$ such that one can take $\tau_B=\kappa_{B,A}\max\{\rho_B/\bar\delta,1\}$.
  \item\label{lem:eb_affine} Each $g_i$ is affine. In this case, $\tau_B$ can be chosen independently of $B$.
\end{enumerate}
\end{lemma}

The proof of Lemma~\ref{lem:errorbound_suff} is given in the Appendix. Case~\ref{lem:eb_affine}) covers coupling equality constraints $Ex=e$, written as $Ex-e\le 0$ and $e-Ex\le 0$, for which Slater's condition always fails. Without any regularity, (A5) may fail: for $n=m=1$, $g_1(x)=x^2$ and $A=[-1,1]$, we have $S=\{0\}$ and $\dist(x,S)/\max\{[g_1(x)]_+,\dist(x,A)\}=1/|x|$ for $0<|x|\le 1$, which is unbounded near $0$.

\subsection{Variational equilibria}\label{subsec:vgne}

By (A2) $S$ is nonempty, convex and compact. A \emph{variational GNE} (v-GNE) is a solution of $\mathrm{VI}(S,F)$, and we let $Z^\star := \SOL(S,F)$. Every v-GNE is a GNE, and, under a constraint qualification such as Slater's condition~\eqref{eq:slater}, v-GNEs are exactly the GNEs at which the KKT multipliers of the coupling constraints coincide for all players \cite{facchinei2007vi}. 

\begin{prop}[{\cite[Secs.~2.2--2.3]{FaccPang1}}]\label{prop:existence}
Under (A1)--(A3), $Z^\star$ is nonempty, convex and compact. Hence, the minimum-norm v-GNE
\begin{equation}\label{eq:min_norm}
  x^\dagger := \Pi_{Z^\star}(0) = \operatorname*{argmin}_{x\in Z^\star}\|x\|
\end{equation}
exists and is unique.
\end{prop}

If $F$ is strongly monotone, $Z^\star$ is a singleton; under (A1), it may be a continuum, and the limit of an iterative method generally depends on its initialization. The minimum-norm v-GNE provides a selection criterion independent of both the algorithm and its initialization. 

In this work, we are interested in the computation of a minimum-norm v-GNE under Assumption~\ref{assum:standing}. 

\subsection{Problem statement}\label{subsec:problem}

Our goal is to compute the minimum-norm v-GNE $x^\dagger$ in~\eqref{eq:min_norm}. Since the v-GNEs are the solutions of $\mathrm{VI}(S,F)$, this amounts to computing the minimum-norm element of $\SOL(S,F)$, a monotone VI over the convex set $S$ defined by the possibly nonlinear coupling constraints. Classical projection methods for $\mathrm{VI}(S,F)$ \cite[Ch.~12]{FaccPang2} are unattractive here: for nonlinear $g_i$, $\Pi_S$ requires the solution of a nonlinear convex program at every iteration. Alternatively, the coupling constraints can be dualized by lifting $\mathrm{VI}(S,F)$ to the primal-dual variable $(x,\lambda)\in A\times\R^m_+$, with operator $(x,\lambda)\mapsto\big(F(x)+\sum_{i\in\mathcal{M}}\lambda_i\nabla g_i(x),\,-g(x)\big)$, as in the extended game with a dual player \cite{TatKamECC25}. This operator is monotone, but, for nonlinear $g_i$, the product of $\lambda$ and $\nabla g_i(x)$ makes it not Lipschitz continuous on the unbounded set $A\times\R^m_+$. Thus, constant step sizes are not justified, and line-search or adaptive step-size techniques, such as the golden ratio method \cite{malitsky2020golden}, are required. Moreover, lifting increases the problem dimension from $n$ to $n+m$, which may slow down convergence, as illustrated in Section~\ref{sec:numerics}.

In view of this discussion, we aim at a method that stays in the primal space, i.e., introduces no dual variables, and avoids projections onto the convex set $S$ of complex structure. To this end, we use only projections onto simple sets: $\Pi_A$ is explicit by~\eqref{eq:proj_A}, and the projection onto the intersection of $A$ with a halfspace is a continuous quadratic knapsack problem, solvable exactly in $\mathcal{O}(n\log n)$ operations \cite{kiwiel2008breakpoint}. If the halfspace contains $X_i$, this intersection is a polyhedral outer approximation of $S$. Projections onto such sets enforce the coupling constraints only asymptotically, and (A5) quantifies how small violations translate into proximity to $S$. This motivates the following problem.

\begin{problem}\label{prob:main}
Under (A1)--(A5), design a (possibly randomized) algorithm that generates random vectors $(x^k)_{k\in\N}$ such that
\begin{enumerate}
  \item\label{prob:conv} $\lim_{k\to\infty}\E\big[\|x^k-x^\dagger\|^2\big]=0$;
  \item\label{prob:local} $x^k\in A$ for every $k\ge 1$;
  \item\label{prob:cost} every iteration $k$ requires only evaluations of $F$, of one coupling function $g_{i_k}$ and of its gradient, and projections onto $A$ intersected with a halfspace containing $X_{i_k}$ (Section~\ref{subsec:linearization}).
\end{enumerate}
\end{problem}

Since $x^\dagger\in S$, item~\ref{prob:conv} implies $\E[\dist(x^k,S)]\to 0$, although $x^k$ need not belong to $X$ for any finite $k$.

\subsection{Linearization of the coupling constraints}\label{subsec:linearization}

For nonlinear $g_i$, even the projection onto a single set $A\cap X_i$ has no closed form. To meet item~\ref{prob:cost}) of Problem~\ref{prob:main}, we therefore replace $X_i$ by the halfspace obtained by linearizing $g_i$ at the current point, as in relaxed projection methods \cite{fukushima1986relaxed}. By convexity, this halfspace contains $X_i$, so that the projection onto it does not cut off any point of $S$.

\begin{lemma}[Linearization]\label{lem:linearization}
Let (A2) and (A4) hold and, for $x,y\in\R^n$ and $i\in\mathcal{M}$, let
\begin{equation}\label{eq:linearization}
  H_i(y) := \{x\in\R^n : g_i(y)+\inner{\nabla g_i(y)}{x-y}\le 0\}.
\end{equation}
\begin{enumerate}
  \item\label{lem:lin_outer} $X_i\subseteq H_i(y)$, and hence $S\subseteq A\cap H_i(y)$.
  \item\label{lem:lin_grad} If $g_i(y)>0$, then $\nabla g_i(y)\neq 0$.
\end{enumerate}
\end{lemma}

\begin{proof}
\ref{lem:lin_outer}) follows from the gradient inequality $g_i(x)\ge g_i(y)+\inner{\nabla g_i(y)}{x-y}$. If $g_i(y)>0$ and $\nabla g_i(y)=0$, the same inequality gives $g_i(x)\ge g_i(y)>0$ for all $x$, i.e., $X_i=\emptyset$, which contradicts (A2); this proves~\ref{lem:lin_grad}).
\end{proof}

By Lemma~\ref{lem:linearization}, if $g_i(y)>0$, then $H_i(y)$ is a halfspace not containing $y$, and
\begin{equation}\label{eq:halfspace_proj}
\begin{gathered}
    \Pi_{H_i(y)}(y) = y-\frac{g_i(y)}{\|\nabla g_i(y)\|^2}\nabla g_i(y),\\
    \dist\big(y,H_i(y)\big) = \frac{g_i(y)}{\|\nabla g_i(y)\|},
\end{gathered}
\end{equation}
whereas $y\in H_i(y)$ if $g_i(y)\le 0$. Moreover, for every bounded set $B$,
\begin{equation}\label{eq:upperG}
    G_B := \sup_{x\in B}\max_{i\in\mathcal{M}}\|\nabla g_i(x)\| < \infty,
\end{equation}
and therefore $\dist(y,H_i(y))\ge[g_i(y)]_+/G_B$ for all $y\in B$.

\section{Main results}\label{sec:main}

\subsection{The algorithm}\label{subsec:algorithm}

For $\varepsilon>0$, let
\begin{equation}\label{eq:F_eps}
  F_\varepsilon := F+\varepsilon I,
\end{equation}
the pseudo-gradient of the game with costs $J_\nu+\frac{\varepsilon}{2}x_\nu^2$. By (A1), $F_\varepsilon$ is strongly monotone on $A$ with modulus $\varepsilon$. Our method, summarized in Algorithm~\ref{alg:main}, alternates a forward step on $F_{\varepsilon_k}$ with a projection onto $A$ intersected with the linearization of one randomly selected coupling constraint.

\begin{algorithm}[!ht]
\caption{Regularized random linearized projection}\label{alg:main}
\begin{algorithmic}[1]
\Require $x^0\in A$; positive sequences $(\alpha_k)_{k\in\N}$ and $(\varepsilon_k)_{k\in\N}$
\For{$k=0,1,2,\ldots$}
  \State draw $i_k\in\mathcal{M}$ uniformly at random, independently of $i_0,\ldots,i_{k-1}$
  \State\label{alg:forward} $y^k \gets x^k-\alpha_k F_{\varepsilon_k}(x^k)$, with $F_{\varepsilon_k}$ as in~\eqref{eq:F_eps}
  \State\label{alg:halfspace} evaluate $g_{i_k}(y^k)$ and $\nabla g_{i_k}(y^k)$, and form the halfspace $H_{i_k}(y^k)$ in~\eqref{eq:linearization}
  \State\label{alg:set} $C_k \gets A\cap H_{i_k}(y^k)$, with the box $A$ of the local constraints in~\eqref{eq:proj_A}
  \State\label{alg:proj} $x^{k+1} \gets \Pi_{C_k}(y^k)$
\EndFor
\end{algorithmic}
\end{algorithm}

Algorithm~\ref{alg:main} differs from the Tikhonov-regularized projection method, $x^{k+1}=\Pi_S(y^k)$, only in that $S$ is replaced by $C_k$. Since $S\subseteq C_k\subseteq A$ by Lemma~\ref{lem:linearization}, $C_k$ is nonempty, closed and convex, the projection onto $C_k$ does not increase the distance to any point of $S$ by~\eqref{eq:fejer_basic}, and every iterate satisfies the local constraints, but not necessarily the coupling ones. Each iteration evaluates $F$ once and $g_{i_k}$ and $\nabla g_{i_k}$ once, so that its cost does not depend on $m$, as in incremental constraint projection methods \cite{wang2015incremental}. If $\Pi_A(y^k)\in H_{i_k}(y^k)$, then $x^{k+1}=\Pi_A(y^k)$; otherwise, $x^{k+1}$ solves a continuous quadratic knapsack problem \cite{kiwiel2008breakpoint}. A fully explicit variant is given in Remark~\ref{rem:sequential}.

Let $\mathcal{F}_k := \sigma(i_0,\ldots,i_{k-1})$ ($\mathcal{F}_0$ trivial) and $\E_k[\,\cdot\,] := \E[\,\cdot\mid\mathcal{F}_k]$. Then $x^k$ and $y^k$ are $\mathcal{F}_k$-measurable and take finitely many values, and, since $i_k$ is uniform and independent of $\mathcal{F}_k$,
\begin{equation}\label{eq:cond_exp}
    \E_k\big[\dist\big(y^k,H_{i_k}(y^k)\big)^2\big] = \frac{1}{m}\sum_{i\in\mathcal{M}}\dist\big(y^k,H_i(y^k)\big)^2\quad\text{a.s.}
\end{equation}
All results below remain valid if $i_k$ is drawn from a distribution $(p_i)_{i\in\mathcal{M}}$ with $p_i>0$.

\subsection{Auxiliary results}\label{subsec:auxiliary}

Let $(\alpha_k)_{k\in\N}$ and $(\varepsilon_k)_{k\in\N}$ be bounded, with $\bar\alpha := \sup_k\alpha_k$ and $\bar\varepsilon := \sup_k\varepsilon_k$, and let
\begin{equation}\label{eq:B}
\begin{gathered}
    B := \{y\in\R^n : \dist(y,A)\le\bar\alpha M\}, \\
    M := \max_{x\in A}\|F(x)\|+\bar\varepsilon\max_{x\in A}\|x\|,
\end{gathered}
\end{equation}
which are, respectively, compact and finite by (A2) and (A3). Since $x^0\in A$ and $x^{k+1}\in C_k\subseteq A$, for every realization of $(i_k)_{k\in\N}$ and every $k\in\N$,
\begin{equation}\label{eq:bounded}
    x^k\in A,\ \ \|F_{\varepsilon_k}(x^k)\|\le M,\ \ \|y^k-x^k\|\le\alpha_k M,\ \ y^k\in B.
\end{equation}
We fix $G>0$ with $G\ge G_B$ in~\eqref{eq:upperG} and the constant $\tau_B$ of (A5) for this set $B$. We further need two known results.

\begin{lemma}[{\cite[Lemma~2.5]{xu2002iterative}}]\label{lem:recursion}
Let $u_k\ge 0$ satisfy $u_{k+1}\le(1-\tau_k)u_k+\tau_k\zeta_k+c_k$ for all $k\ge k_0$, where $\tau_k\in(0,1)$, $\sum_k\tau_k=\infty$, $\limsup_{k}\zeta_k\le 0$, $c_k\ge 0$ and $\sum_k c_k<\infty$. Then $u_k\to 0$.
\end{lemma}

\begin{lemma}[Tikhonov path]\label{lem:TikSeq2}
Let (A1)--(A3) hold. For every $\varepsilon>0$, $\mathrm{VI}(S,F_\varepsilon)$ has a unique solution $x(\varepsilon)$, and
\begin{enumerate}
  \item\label{lem:tik_bound} $\|x(\varepsilon)\|\le\|x^\dagger\|$ for every $\varepsilon>0$;
  \item\label{lem:tik_limit} $\lim_{\varepsilon\downarrow 0}x(\varepsilon)=x^\dagger$;
  \item\label{lem:tik_lip} $\|x(\varepsilon)-x(\varepsilon')\| \le\|x^\dagger\|\,|\varepsilon-\varepsilon'|/\varepsilon$ for all $\varepsilon,\varepsilon'>0$.
\end{enumerate}
\end{lemma}

Items~\ref{lem:tik_bound}) and~\ref{lem:tik_limit}) are classical \cite{BakPolyak74}, \cite[Sec.~12.2]{FaccPang2}. Item~\ref{lem:tik_lip}) is obtained as in \cite{kannan2012distributed,koshal2013regularized}, by adding the variational inequalities defining $x(\varepsilon)$ and $x(\varepsilon')$ and bounding $\|x(\varepsilon')\|$ by~\ref{lem:tik_bound}).

\subsection{Asymptotic feasibility}\label{subsec:feasibility}

The iterates of Algorithm~\ref{alg:main} need not satisfy the coupling constraints, but, as the next lemma shows, their infeasibility is of the order of the step size: in expectation, the random projections reduce the distance to $S$ by a fixed fraction, whereas each forward step increases it by at most $\alpha_k M$.

\begin{lemma}[Asymptotic feasibility]\label{lem:feasibility}
Let (A2)--(A5) hold, and let $(\alpha_k)_{k\in\N}$ and $(\varepsilon_k)_{k\in\N}$ be positive and bounded, with
\begin{equation}\label{eq:ratio}
  \lim_{k\to\infty}\alpha_{k-1}/\alpha_k=1.
\end{equation}
Then there is $D>0$ such that, for all $k\in\N$,
\begin{equation}\label{eq:feasibility}
  \E\big[\dist(x^k,S)^2\big]\le D^2\alpha_k^2, \quad \E\big[\dist(x^k,S)\big]\le D\,\alpha_k.
\end{equation}
\end{lemma}

\begin{proof}
\emph{Step~1.} Since $\Pi_S(y^k)\in S\subseteq C_k$, \eqref{eq:fejer_basic} gives
\begin{equation}\label{eq:fejer}
  \dist(x^{k+1},S)^2 \le\dist(y^k,S)^2-\dist(y^k,C_k)^2.
\end{equation}

\emph{Step~2.} We bound the decrease $\dist(y^k,C_k)^2$ in~\eqref{eq:fejer} from below, in conditional expectation. Fix $i\in\mathcal{M}$. If $g_i(y^k)\le 0$, then $y^k\in H_i(y^k)$. Otherwise, \eqref{eq:halfspace_proj} and $\|\nabla g_i(y^k)\|\le G$, which holds since $y^k\in B$ by~\eqref{eq:bounded}, give $\dist(y^k,H_i(y^k))=g_i(y^k)/\|\nabla g_i(y^k)\|\ge g_i(y^k)/G$. In both cases, $\dist(y^k,H_i(y^k))\ge[g_i(y^k)]_+/G$. Hence, by~\eqref{eq:cond_exp}, and since a sum of nonnegative terms is not smaller than its largest term, a.s.,
\begin{align*}
  \E_k\big[\dist(y^k,H_{i_k}(y^k))^2\big]
  \ge\frac{\max_{i}[g_i(y^k)]_+^2}{mG^2}.
\end{align*}
Since $C_k=A\cap H_{i_k}(y^k)$ is contained both in $A$ and in $H_{i_k}(y^k)$, we have $\dist(y^k,C_k)\ge\max\{\dist(y^k,A),\dist(y^k,H_{i_k}(y^k))\}$. Moreover, $\dist(y^k,A)$ is $\mathcal{F}_k$-measurable, because $y^k$ is. Taking $\E_k$ and using the last display and (A5) at $y^k\in B$, we obtain, a.s.,
\begin{align*}
  \E_k\big[\dist(y^k,C_k)^2\big]
  &\ge\max\Big\{\dist(y^k,A)^2,\ \frac{\max_{i}[g_i(y^k)]_+^2}{mG^2}\Big\}\\
  &\ge c\,\dist(y^k,S)^2,
  \end{align*}
where $c := \tau_B^{-2}\min\big\{1,\ (mG^2)^{-1}\big\}\in(0,1]$ since $\tau_B\ge 1$. Finally, $\dist(y^k,S)$ is $\mathcal{F}_k$-measurable, so that taking $\E_k$ in~\eqref{eq:fejer} gives, a.s.,
\begin{equation}\label{eq:contraction}
  \E_k\big[\dist(x^{k+1},S)^2\big]\le q\,\dist(y^k,S)^2,\quad q := 1-c\in[0,1).
\end{equation}
In words, in conditional expectation, each projection removes at least the fraction $c$ of the squared distance of $y^k$ to $S$.

\emph{Step~3.} We now account for the forward step, which may increase the distance to $S$. By line~\ref{alg:forward} of Algorithm~\ref{alg:main} and~\eqref{eq:bounded}, $\|y^k-x^k\|=\alpha_k\|F_{\varepsilon_k}(x^k)\|\le\alpha_k M$. Since $\dist(\cdot,S)$ is $1$-Lipschitz continuous, $\dist(y^k,S)\le\dist(x^k,S)+\alpha_k M$. For every $\theta>0$, Young's inequality $2ab\le\theta a^2+\theta^{-1}b^2$ gives $(a+b)^2\le(1+\theta)a^2+(1+\theta^{-1})b^2$. Applying it with $\theta := (1-q)/2>0$ in~\eqref{eq:contraction}, we obtain, a.s.,
by taking total expectations, that $u_k := \E[\dist(x^k,S)^2]$ satisfies
\begin{equation}\label{eq:recursion2}
  u_{k+1}\le\bar q\,u_k+\Gamma\alpha_k^2,
\end{equation}
where $\Gamma := q(1+\theta^{-1})M^2$ and $\bar q := q(1+\theta)=q(3-q)/2$. Since $\bar q-1=(q-1)(2-q)/2$ and $q\in[0,1)$, we have $\bar q<1$. Thus, $\theta$ is chosen small enough to preserve the contraction of Step~2 after the perturbation by the forward step.

\emph{Step~4.} For $k\ge 1$, let $s_k := \alpha_{k-1}/\alpha_k$ and $v_k := u_k/\alpha_{k-1}^2$. Dividing~\eqref{eq:recursion2} by $\alpha_k^2$ gives $v_{k+1}\le\bar q s_k^2 v_k+\Gamma$. By~\eqref{eq:ratio}, there is $k_1$ with $\bar q s_k^2\le\hat q := (1+\bar q)/2<1$ for $k\ge k_1$, so that, by induction, $v_k\le V := \max\{v_{k_1},\Gamma/(1-\hat q)\}$ and $u_k=s_k^2\alpha_k^2v_k\le\bar s^2V\alpha_k^2$ for $k\ge k_1$, where $\bar s := \sup_k s_k<\infty$. The claim follows with $D^2 := \max\{\bar s^2 V,\max_{k<k_1}u_k/\alpha_k^2,1\}$ and Jensen's inequality.
\end{proof}

Lemma~\ref{lem:feasibility} establishes the asymptotic feasibility of the iterates of Algorithm~\ref{alg:main}: although $x^k$ need not satisfy the coupling constraints, its expected distance to $S$ vanishes at the rate $\mathcal{O}(\alpha_k)$. This result is required in the convergence analysis of the next subsection, where the infeasibility of $x^k$ enters as an error term that must be controlled by the step size.

\subsection{Convergence to the minimum-norm v-GNE}\label{subsec:convergence}

The iterates chase a moving target: each step brings $x^{k+1}$ closer to $x(\varepsilon_k)$, which converges to $x^\dagger$. The parameters must be balanced so that the iterates keep up with the target, while the errors due to the forward step and to the infeasibility remain summable.

\begin{theorem}\label{thm:main}
Let (A1)--(A5) hold, and let $(\alpha_k)_{k\in\N}$ and $(\varepsilon_k)_{k\in\N}$ be positive sequences such that
\begin{subequations}\label{eq:step}
\begin{align}
  \lim_{k\to\infty}\frac{\alpha_{k-1}}{\alpha_k}=1,
  \qquad \sum_{k=0}^{\infty}\alpha_k^2<\infty,
  \label{eq:step_alpha}\\
  \lim_{k\to\infty}\varepsilon_k=0,
  \qquad \sum_{k=0}^{\infty}\alpha_k\varepsilon_k=\infty,
  \\\lim_{k\to\infty}\frac{|\varepsilon_{k-1}-\varepsilon_k|}{\alpha_k\varepsilon_k^2}=0.
  \label{eq:step_eps}
\end{align}
\end{subequations}
Then the iterates of Algorithm~\ref{alg:main} satisfy $\lim_{k\to\infty}\E\big[\|x^k-x^\dagger\|^2\big]=0$ and $\E[\dist(x^k,S)]=\mathcal{O}(\alpha_k)$; in particular, $x^k\to x^\dagger$ in probability.
\end{theorem}

\begin{proof}
By~\eqref{eq:step}, $\alpha_k\to 0$ and $\varepsilon_k\to 0$, so~\eqref{eq:bounded} and Lemma~\ref{lem:feasibility} apply. Let $\tau_k := \alpha_k\varepsilon_k$, $\Delta_k := |\varepsilon_{k-1}-\varepsilon_k|$, and let $k_0\ge 1$ be such that $\tau_k\le\frac12$ for $k\ge k_0$.

\emph{Step~1.} Let $z := x(\varepsilon_k)$. Since $z\in S\subseteq C_k$, \eqref{eq:fejer_basic} and lines~\ref{alg:forward} and~\ref{alg:proj} of Algorithm~\ref{alg:main} give
\begin{align}\label{eq:dist1}
  \|x^{k+1}-z\|^2\le\|x^k-z\|^2-2\alpha_k\inner{F_{\varepsilon_k}(x^k)}{x^k-z} \nonumber\\
   +\alpha_k^2\|F_{\varepsilon_k}(x^k)\|^2.
\end{align}
We write $F_{\varepsilon_k}(x^k)=[F(x^k)-F(z)]+\varepsilon_k(x^k-z)+F_{\varepsilon_k}(z)$. First, $\inner{F(x^k)-F(z)}{x^k-z}\ge 0$ by (A1), since $x^k,z\in A$; here monotonicity on $A$, not only on $S$, is needed. Second, $\inner{\varepsilon_k(x^k-z)}{x^k-z}=\varepsilon_k\|x^k-z\|^2$. Third, since $z$ solves $\mathrm{VI}(S,F_{\varepsilon_k})$, $\Pi_S(x^k)\in S$ and $\|F_{\varepsilon_k}(z)\|\le M$,
\begin{align*}
  \inner{F_{\varepsilon_k}(z)}{x^k-z}
  &\ge\inner{F_{\varepsilon_k}(z)}{x^k-\Pi_S(x^k)}\\
  &\ge -M\dist(x^k,S).
\end{align*}
Inserting these bounds and~\eqref{eq:bounded} into~\eqref{eq:dist1}, taking expectations, and using Lemma~\ref{lem:feasibility}, we obtain
\begin{equation}\label{eq:dist2}
  \E\|x^{k+1}-x(\varepsilon_k)\|^2
  \le(1-2\tau_k)\E\|x^k-x(\varepsilon_k)\|^2+\Lambda\alpha_k^2,
\end{equation}
where $\Lambda := M^2+2MD$.

\emph{Step~2.} Let $u_k := \E\|x^k-x(\varepsilon_{k-1})\|^2$. For $k\ge k_0$, $(a+b)^2\le(1+\tau_k)a^2+(1+\tau_k^{-1})b^2$ and Lemma~\ref{lem:TikSeq2}(\ref{lem:tik_lip}) give
\begin{align*}
  \|x^k-x(\varepsilon_k)\|^2
  &\le (1+\tau_k)\|x^k-x(\varepsilon_{k-1})\|^2\\
  &\quad+(1+\tau_k^{-1})\|x^\dagger\|^2\Delta_k^2/\varepsilon_k^2.
\end{align*}
Inserting this into~\eqref{eq:dist2} and using $(1-2\tau_k)(1+\tau_k)\le 1-\tau_k$ and $(1-2\tau_k)(1+\tau_k^{-1})\le\tau_k^{-1}$, we obtain, for $k\ge k_0$,
\begin{equation}\label{eq:2f}
  u_{k+1}\le(1-\tau_k)u_k+\tau_k\zeta_k+\Lambda\alpha_k^2,
  \ \
  \zeta_k := \|x^\dagger\|^2\Big(\frac{\Delta_k}{\alpha_k\varepsilon_k^2}\Big)^2.
\end{equation}

\emph{Step~3.} By~\eqref{eq:step}, $\sum_k\tau_k=\infty$, $\zeta_k\to 0$ and $\sum_k\alpha_k^2<\infty$, so that $u_k\to 0$ by Lemma~\ref{lem:recursion}. Finally, $\E\|x^k-x^\dagger\|^2\le 2u_k+2\|x(\varepsilon_{k-1})-x^\dagger\|^2\to 0$ by Lemma~\ref{lem:TikSeq2}(\ref{lem:tik_limit}). The remaining claims follow from Lemma~\ref{lem:feasibility} and Markov's inequality.
\end{proof}

Algorithm~\ref{alg:main} solves Problem~\ref{prob:main}: item~\ref{prob:conv} holds by Theorem~\ref{thm:main}, items~\ref{prob:local} and~\ref{prob:cost} by construction. Indeed, $x^{k+1}=\Pi_{C_k}(y^k)\in C_k\subseteq A$ by lines~\ref{alg:set} and~\ref{alg:proj}, and each iteration evaluates $F$, $g_{i_k}$ and $\nabla g_{i_k}$ once and projects onto $C_k=A\cap H_{i_k}(y^k)$, where $H_{i_k}(y^k)\supseteq X_{i_k}$ by Lemma~\ref{lem:linearization}.

\begin{remark}[Parameters]\label{rem:step}
The sequences $\alpha_k=\alpha_0(k+1)^{-a}$ and $\varepsilon_k=\varepsilon_0(k+1)^{-b}$, $\alpha_0,\varepsilon_0>0$, satisfy~\eqref{eq:step} whenever $\frac12<a<1$ and $0<b<1-a$, since then $\Delta_k/(\alpha_k\varepsilon_k^2)=\mathcal{O}(k^{a+b-1})$. Theorem~\ref{thm:main} provides no rate, since no estimate of $\|x(\varepsilon)-x^\dagger\|$ is available under (A1)--(A5); the slowly vanishing $\varepsilon_k$ is the price of selecting $x^\dagger$. Note that Lipschitz continuity of $F$ is not required: $F$ enters the analysis only through (A1) and the bound $M$ in~\eqref{eq:B}, which is finite due to (A3) since $A$ is compact, and the error $\alpha_k^2M^2$ of the forward step is summable by~\eqref{eq:step_alpha}. Hence, unlike the primal--dual methods discussed in Section~\ref{subsec:problem}, Algorithm~\ref{alg:main} needs neither a Lipschitz constant nor line-search or adaptive step-size rules.
\end{remark}

\begin{remark}[A fully explicit variant]\label{rem:sequential}
The projection onto $C_k$ can be replaced by $x^{k+1}=\Pi_A\big(\Pi_{H_{i_k}(y^k)}(y^k)\big)$, which is explicit by~\eqref{eq:proj_A} and~\eqref{eq:halfspace_proj} and, except for the scalar $g_{i_k}(y^k)/\|\nabla g_{i_k}(y^k)\|^2$, componentwise, which is convenient for distributed implementations. Lemma~\ref{lem:feasibility} and Theorem~\ref{thm:main} remain valid with $c$ replaced by $c/2$: both projections are nonexpansive and fix every point of $S$, and applying~\eqref{eq:fejer_basic} twice and using that $\dist(\cdot,A)$ is $1$-Lipschitz, the decrease in~\eqref{eq:fejer} is at least $\frac12\max\{\dist(y^k,H_{i_k}(y^k)),\dist(y^k,A)\}^2$.
\end{remark}

\begin{remark}[Nonsmooth coupling constraints]\label{rem:nonsmooth}
Since each $g_i$ is convex and finite, $\partial g_i(y)$ is nonempty, and bounded uniformly on bounded sets \cite[Thms.~23.4, 24.7]{rockafellar1970convex}. Lemma~\ref{lem:linearization} and~\eqref{eq:halfspace_proj} use only the subgradient inequality, and (A5) and Lemma~\ref{lem:errorbound_suff} involve no derivatives. Hence, all results hold without (A4), with $\nabla g_i(y)$ replaced by any subgradient.
\end{remark}

\section{Numerical examples}\label{sec:numerics}

We compare Algorithm~\ref{alg:main} with primal--dual formulations of the v-GNE problem, to show its advantages in the two respects for which it was designed: the selection of $x^\dagger$ (Section~\ref{subsec:num_selection}) and that the cost per iteration does not grow with the number $m$ of coupling constraints (Section~\ref{subsec:num_cost}). As primal--dual methods, we use the state-of-the-art adaptive methods aGRAAL \cite{malitsky2020golden} and PF-NE-EG \cite{shen2026parameterfree}, applied to the lifted operator of Section~\ref{subsec:problem}, $\mathcal{B}(x,\lambda):=\big(F(x)+\sum_{i\in\mathcal{M}}\lambda_i\nabla g_i(x),\, -g(x)\big)$ on $A\times\R^m_+$. Both converge to v-GNEs, and every evaluation of $\mathcal{B}$ involves all $m$ constraints. The experiments illustrate the asymptotic analysis of Section~\ref{sec:main}.

\subsection{Experimental protocol}\label{subsec:num_setup}

We used MATLAB R2025b on an Intel Core Ultra~7 processor with 32\,GB of memory (Windows~11 Pro). Algorithm~\ref{alg:main} uses the exact breakpoint projection of Section~\ref{subsec:algorithm}, Mersenne-Twister index sequences, and $\alpha_k=\alpha_0(k+1)^{-a}$, $\varepsilon_k=\varepsilon_0(k+1)^{-b}$. The exponents $(a,b)$, which satisfy~\eqref{eq:step} by Remark~\ref{rem:step}, and the constants $\alpha_0,\varepsilon_0$ are chosen for each example and stated there. aGRAAL ($\varphi=1.5$) and PF-NE-EG ($\theta=\rho=0.9$) follow \cite[Alg.~1]{malitsky2020golden} and \cite[Alg.~2]{shen2026parameterfree}, start from $(x^0,0)$, and use a dual scaling $\lambda=\kappa\nu$. We count one evaluation of $F$ per iteration of Algorithm~\ref{alg:main} and per evaluation of $\mathcal{B}$. The stopping criteria are stated with each example.

\subsection{Equilibrium selection}\label{subsec:num_selection}

\begin{example}[Minimum-loss load allocation]\label{ex:selection}
Two power plants with outputs $x_\nu\in A_\nu=[0,3]$ cover a demand $d=2$ and pay the same imbalance charge, $J_1(x)=J_2(x)=\frac12(x_1+x_2-d)^2$. This means that $F(x)=(x_1+x_2-d)\,(1,1)$, which is monotone, but not strictly monotone. The power plants are coupled by a regional emission limit given as
\[
  g_1(x)=2x_1^2+\tfrac12x_2^2-2\le 0,
\] 
where plant~$1$ is assumed to be an older or dirtier unit. Assumption~\ref{assum:standing} holds, in particular, Slater's condition~\eqref{eq:slater} holds with $\bar x=0$.
\end{example}

Every split of the demand that respects the emission limit is a v-GNE, so market rules alone do not determine the load allocation. If both plants feed the load through lines of equal resistance, the transmission losses are proportional to $\|x\|^2$. This gives the selection criterion a physical justification: among all equilibria, the minimum-norm v-GNE is the one with the smallest transmission losses, and hence the one a system operator would want to implement \cite[Sec.~3.3]{frank2016introduction}.

In this example, both $Z^\star$ and $x^\dagger$ are explicit. Since $F=\nabla\Phi$ with $\Phi(x)=\frac12(x_1+x_2-d)^2$ convex, the v-GNEs are the minimizers of $\Phi$ over $S$, so that $Z^\star=\{(t,2-t)\mid 0\le t\le\frac45\}$ and $x^\dagger=(\frac45,\frac65)$, at which the emission limit is active. Over $Z^\star$, the losses range from $\|x^\dagger\|^2=2.08$ to $4$ at $(0,2)$, so the selection matters. Moreover, $x(\varepsilon)\neq x^\dagger$ for every $\varepsilon>0$, since $-F_\varepsilon(x^\dagger)=-\varepsilon x^\dagger$ does not belong to the normal cone $\Normal_S(x^\dagger)$, generated by $\nabla g_1(x^\dagger)=(\frac{16}{5},\frac65)$.

To compare the methods, we run aGRAAL and PF-NE-EG with $\kappa=1$ and Algorithm~\ref{alg:main} for $20$ initial points drawn uniformly from $A$ (seed $1$), and measure their accuracy on $Z^\star$ by $\mathrm{res}(x):=|x_1+x_2-d|/\sqrt2+[g_1(x)]_+$, which vanishes exactly on $Z^\star$, for this example. The primal--dual methods stop when $\mathrm{res}\le10^{-4}$, whereas Algorithm~\ref{alg:main}, with $(a,b)=(0.51,0.48)$, $\alpha_0=8$ and $\varepsilon_0=0.025$, runs a fixed budget of $K=1.5\cdot10^5$ iterations, chosen to reach a similar accuracy; since $m=1$, it is deterministic here. 

\begin{figure}[t]
  \centering
  \includegraphics[width=\columnwidth]{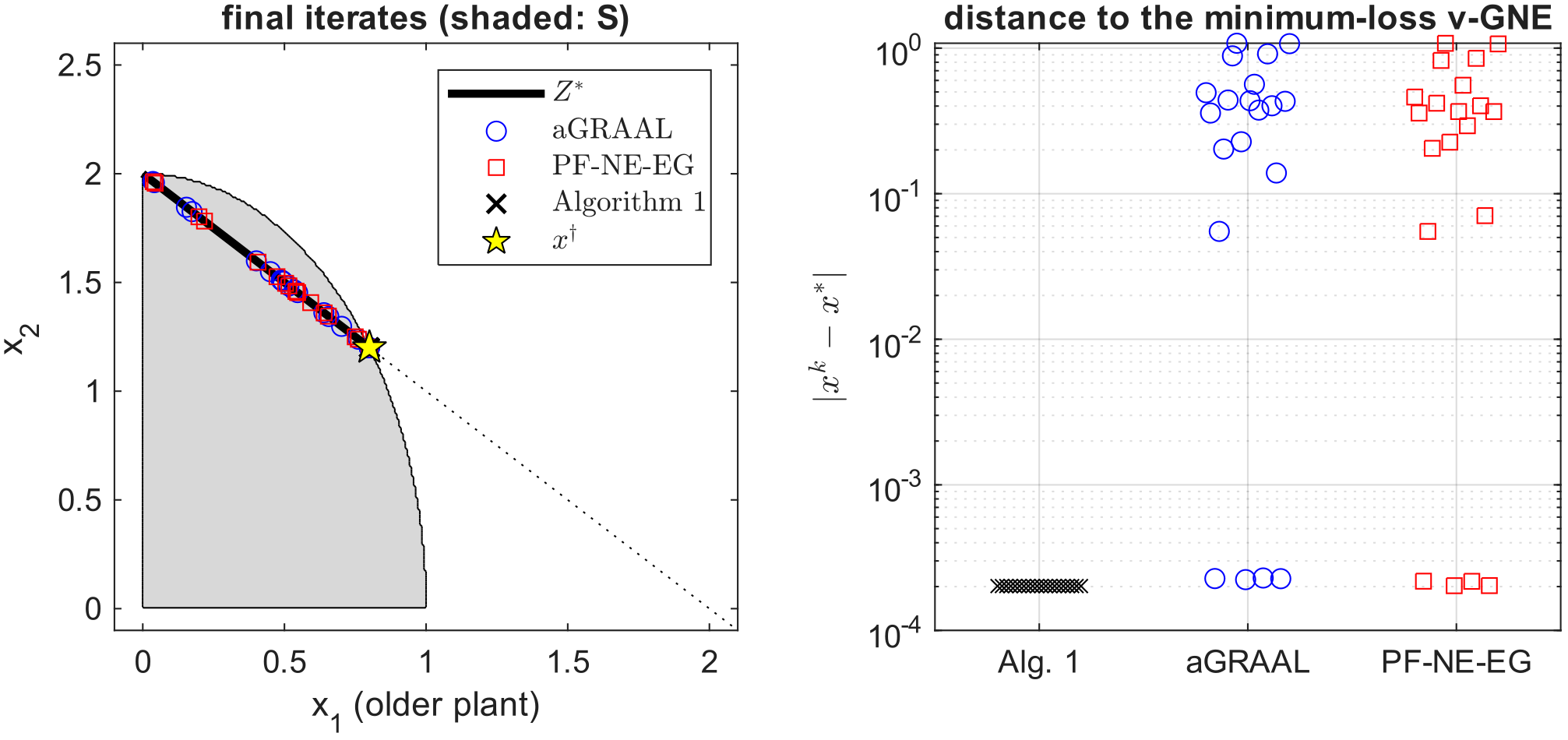}
  \caption{Example~\ref{ex:selection}: final iterates (left; shaded:
  $S$; thick: $Z^\star$; dotted: $x_1+x_2=d$) and their distances to
  $x^\dagger$ (right). Crosses: Algorithm~\ref{alg:main}; circles:
  aGRAAL; squares: PF-NE-EG; star: $x^\dagger$.}
  \label{fig:selection}
\end{figure}

The primal--dual methods stop after a median of $71$ and $61$ evaluations at some point of $Z^\star$ that depends on $x^0$ (Fig.~\ref{fig:selection}). For them, only $4$ of the $20$ runs end near $x^\dagger$. In contrast, Algorithm~\ref{alg:main} approaches $x^\dagger$ from every initial point, thanks to its regularization.

\subsection{Computational cost with many coupling constraints}\label{subsec:num_cost}

\begin{example}[Service market]\label{ex:market} 
Consider $n=30$ suppliers with outputs $x_\nu\in A_\nu=[0,U]$, each serving one or two of $r=6$ markets with supply $Bx$, $B\in\R^{r\times n}_{\ge 0}$, and demand $\tau$, and with costs $J_\nu(x)=\tfrac12 M_{\nu\nu}x_\nu^2+x_\nu\big(\sum_{\mu\neq\nu}M_{\nu\mu}x_\mu-q_\nu\big)$, where $M=B^\top(I+K)B$, $q=B^\top(I+K)\tau$ and $K=-K^\top\neq 0$. Then we have that $F(x)=B^\top(I+K)(Bx-\tau)$, which is monotone, but not strictly (because $B^\top B$ has rank $r$). There are $m=200$ coupling constraints, one per service region, which require that the entropic risk of the unmet demand of the region to be nonpositive, that is,
\[
  g_i(x)=\tfrac{1}{\beta}\log\Big(\tfrac1S\textstyle\sum_{s=1}^S e^{\beta(b_{is}-a_{is}^\top x)}\Big)\le 0,
\]
where $S=1000$ is the number of demand scenarios, $b_{is}$ is the demand of region~$i$ in scenario~$s$, $a_{is}^\top x$ is the service provided by the $3$--$4$ suppliers of the region in that scenario, so that $b_{is}-a_{is}^\top x$ is its unmet demand, and $\beta=1$ is the risk aversion. 

The matrices $B$ and $K$, the vectors $a_{is}$ and the scalars $b_{is}$ are random (seed $1$), except in four regions whose $a_{is}$, $b_{is}$ are constructed, together with $\tau$ and $U$, so that $x^\dagger$ is known in closed form and these four constraints are active at it.
\end{example}

Since $F$ vanishes exactly on $\{Bx=\tau\}$ and, by monotonicity, every v-GNE $x$ satisfies $\|B(x-x^\dagger)\|^2\le 0$, we then have that $Z^\star=\{x\in S\mid Bx=\tau\}$. 
Assumption~\ref{assum:standing} holds, in particular, Slater's condition~\eqref{eq:slater} holds with $\bar x=U\mathbf{1}$. All constraints are violated at $x^0=0$, and the minimum-norm point of $\{Bx=\tau\}$ violates $4$ of them, so they must be enforced along the way.

All methods stop when $\mathrm{res}(x):=\|B^+(Bx- \tau)\|+\max_i[g_i(x)]_+\le\delta$, which, for this example, vanishes exactly on $Z^\star$. aGRAAL and PF-NE-EG check $\mathrm{res}$ after every iteration at no extra cost, since $g(x)$ is part of $\mathcal{B}$. Algorithm~\ref{alg:main} needs all $m$ constraints for it, at about the cost of one evaluation of $\mathcal{B}$, and checks it at iterations $10,12,14,17,\ldots$ (ratio $1.2$), at most $20\%$ late; the checks are included in its time ($13\%$ of it). Algorithm~\ref{alg:main} uses $(a,b)=(0.65,0.3)$, $\alpha_0=2/\|M\|$ and $\varepsilon_0=2\cdot10^{-9}$, with at most $3\cdot10^4$ iterations, and we report the median over $5$ index sequences. The small $\varepsilon_0$ prioritizes fast and accurate convergence to $Z^\star$, at the price that the iterates do not approach $x^\dagger$ within the run. The primal--dual methods use at most $2000$ evaluations and the best $\kappa\in\{0.3,1,3,10\}$ ($\kappa=1$ for both), without counting this tuning, and all methods store the scenario data as sparse matrices.

\begin{figure}[t]
  \centering
  \includegraphics[width=\columnwidth]{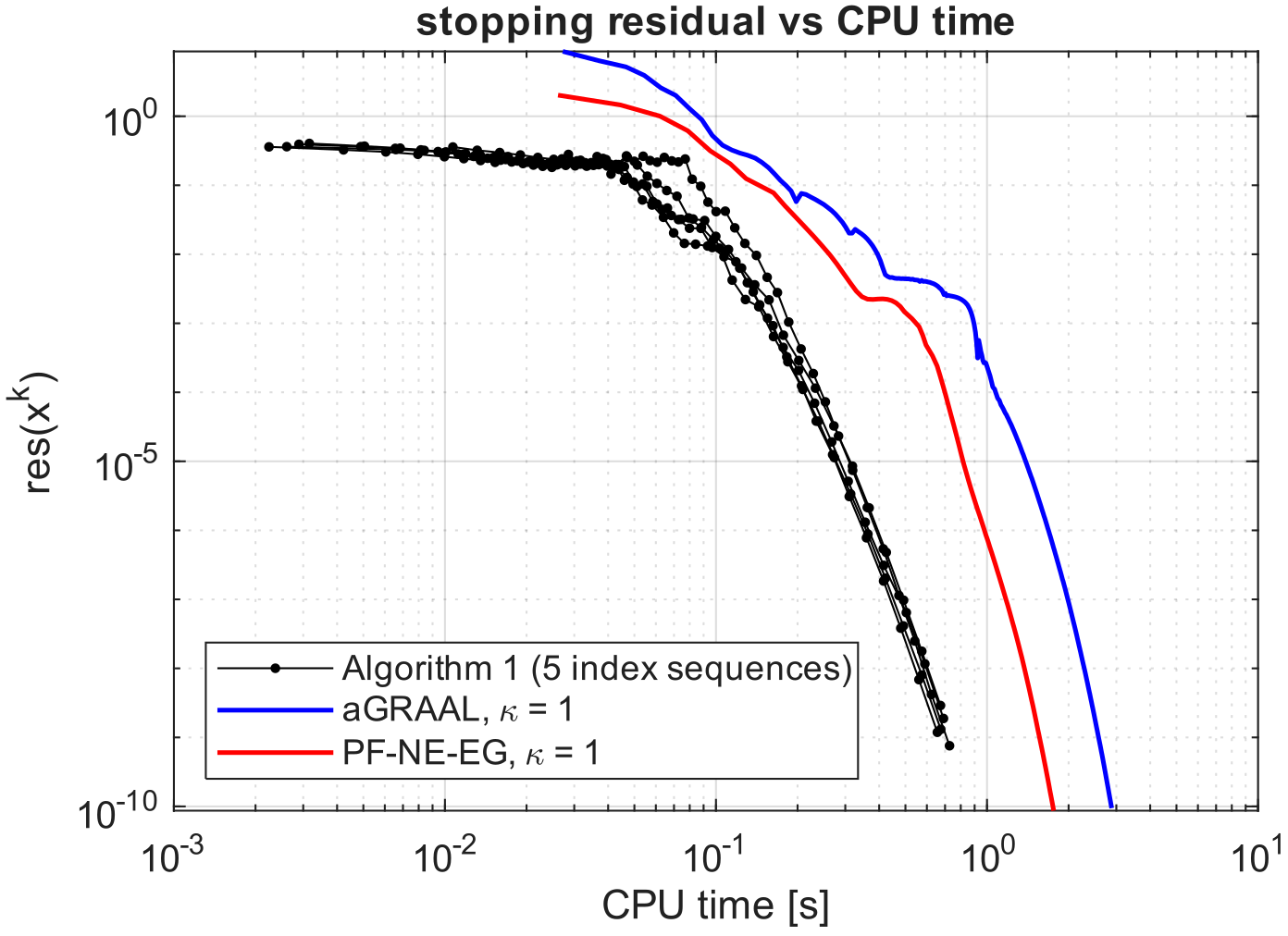}
  \caption{Example~\ref{ex:market}: $\mathrm{res}(x^k)$ versus CPU time;
  Algorithm~\ref{alg:main} at its checks ($5$ index sequences), aGRAAL
  and PF-NE-EG with the best $\kappa$.}
  \label{fig:cpu}
\end{figure}

One iteration of Algorithm~\ref{alg:main} takes $0.022$\,ms and one evaluation of $\mathcal{B}$ about $8.4$\,ms, since it involves all $200$ regions instead of one. Although the primal--dual methods need far fewer steps, Algorithm~\ref{alg:main} stops first at every tolerance and it is $2.0$ to $2.9$ times faster than PF-NE-EG and $3.1$ to $4.5$ times faster than aGRAAL (Fig.~\ref{fig:cpu}).

\section{Conclusion}\label{sec:conclusion}

We proposed a randomized Tikhonov-regularized algorithm for computing the minimum-norm v-GNE of merely monotone games with convex, possibly nonlinear, coupling constraints. Each iteration projects onto the local action sets intersected with the linearization of a single, randomly selected coupling constraint, so that no projection onto the nonlinear feasible set is required. Under an error bound on the constraints, implied, e.g., by Slater's condition, the infeasibility of the iterates is of the order of the step size, and the iterates converge in mean square to the minimum-norm v-GNE. Future work includes convergence rates under additional regularity, non-uniform sampling of the constraints, and fully distributed implementations.

\section*{Declaration of Generative AI Use}

The authors used Claude (Anthropic) to assist with language editing and improving the clarity, grammar, and parallelism of the manuscript text, and to clarify terminology and phrasing during the writing process. The tool was not used to generate original research content, data, analysis, or conclusions. All AI-assisted suggestions were reviewed, verified, and edited by the authors, who take full responsibility for the accuracy, originality, and final content of the manuscript.

\bibliographystyle{IEEEtran}
\bibliography{refs}

\appendix

\section*{Proof of Lemma~\ref{lem:errorbound_suff}}

\ref{lem:eb_slater}) Since $g$ is continuous, $X$ contains a ball centered at $\bar x$; hence, $\bar x\in A\cap\operatorname{int}X$, and $\{X,A\}$ is boundedly linearly regular \cite{BauschkeBorwein1996}, i.e., there is $\kappa_{B,A}\ge 1$ with $\dist(x,S) \le \kappa_{B,A}\max\{\dist(x,X),\dist(x,A)\}$ for all $x\in B$. Moreover, by Robinson's error bound \cite{robinson1975application}, $\dist(x,X) \le (\rho_B/\bar\delta) \max_{i\in\mathcal{M}}[g_i(x)]_+$ for all $x\in B$. Combining both inequalities gives~\eqref{eq:errorbound}.

\ref{lem:eb_affine}) $S$ is a nonempty polyhedron described by the affine inequalities $g(x)\le 0$, $a-x\le 0$ and $x-b\le 0$. By Hoffman's lemma \cite{hoffman}, there is $\tau_H>0$ with $\dist(x,S)\le\tau_H\|([g(x)]_+,[a-x]_+,[x-b]_+)\|$ for all $x\in\R^n$, where $[\cdot]_+$ acts componentwise. By~\eqref{eq:proj_A}, $\|([a-x]_+,[x-b]_+)\|=\dist(x,A)$, and $\|[g(x)]_+\|\le\sqrt{m}\max_{i}[g_i(x)]_+$. Hence,~\eqref{eq:errorbound} holds with $\tau_B=\max\{\sqrt{m+1}\,\tau_H,1\}$ for every $B$.
\hfill$\blacksquare$

\end{document}